\documentclass[11pt]{article}
\usepackage{indentfirst}
\usepackage{graphicx}
\usepackage{amsfonts}
\usepackage{url}
\usepackage{algorithm}
\usepackage{algpseudocode}
\usepackage{listings}
\usepackage{lscape} 
\newcommand{\sn}{\mathrm{sn}}
\newcommand{\cn}{\mathrm{cn}}
\newcommand{\dn}{\mathrm{dn}}

\newtheorem{teor}{Theorem}[section]

\author{E. Scheiber\thanks{e-mail: scheiber@unitbv.ro}}
\title{On Computing Jacobi's Elliptic Function \texttt{sn} }
\date{}
\begin{document}
\maketitle
\begin{abstract}

The paper presents a method to compute the Jacobi's elliptic function \texttt{sn} on the period
parallelogram.  For fixed  $m$ it requires first to compute  the complete elliptic integrals $K=K(m)$ and $K'=K(1-m).$
The Newton method is used to compute $\sn(z,m),$ when $z\in [0,K]\cup[0,i K').$ The computation
in any other point does not require the usage of any numerical procedure, it is done only
with the help of the properties of $\sn.$ 
\\

2010 \textit{Mathematics Subject Classification:} 65D20, 33F05.

\textit{Key words:} elliptic functions, elliptic integrals, arithmetic-geometric mean, 
\end{abstract}

\section{Introduction}
The paper presents a method to compute the Jacobi's elliptic function \texttt{sn} on the period
parallelogram.  For fixed $m\in (0,1)$ it requires first to compute  the complete elliptic integrals $K=K(m)$ and $K'=K(1-m).$
The function to compute the first complete elliptic integral uses the arithmetic-geometric mean, as a
consequence of Gauss's theorem. 

The Newton method to solve a nonlinear algebraic equation is used to compute 
$\sn(z,m),$ when $z\in [0,K]\cup[0,i K').$ The computation
in any other point does not require the usage of any numerical procedure, it is done only
with the help of the properties of $\sn$  and
its values on some points from $[0,K]\cup i[0,K').$ 

The validity of the method is exemplified  with the help of a \textit{Scilab} application. 
The obtained results are very good approximations of the values given by the corresponding
functions from \textit{Scilab} and \textit{Mathematica.}

The computation of the elliptic integrals and of the elliptic functions were studied in many papers, e.g. \cite{3}, \cite{8}, \cite{4},
as well as the included bibliography.

\section{Incomplete elliptic integral of first kind}

The following incomplete and complete elliptic integrals of first kind are defined respectively by, \cite{2},
$$
F(\phi,m)=\int_0^{\sin{\phi}}\frac{\mathrm{d}t}{\sqrt{(1-t^2)(1-mt^2)}}=\int_0^{\phi}\frac{\mathrm{d}\theta}{\sqrt{1-m\sin^2\theta}}
$$
and
$$
K(m)=\int_0^1\frac{\mathrm{d}t}{\sqrt{(1-t^2)(1-mt^2)}}=F(\frac{\pi}{2},m).
$$

In order to compute $F(\phi,m)$ we recall a result established by  Carl Friedrich GAUSS (1777-1855) in 1799, \cite{1}, \cite{5}:

\begin{teor}\label{teor1}
If $a$ and $b$ are positive reals and $M(a,b)$ is their the arithmetic-geometric mean then
\begin{equation}\label{eq9}
\frac{1}{M(a,b)}=\frac{2}{\pi}\int_0^{\frac{\pi}{2}}\frac{\mathrm{d}x}{\sqrt{a^2 \cos^2{x}+b^2\sin^2{x}}}.
\end{equation}
\end{teor}

For $a>b>0$ and $0\le\phi\le\frac{\pi}{2}$ we shall take care of the integral
\begin{equation}\label{eq1}
I(a,b,\phi)=\int_0^{\phi}\frac{\mathrm{d}x}{\sqrt{a^2\cos^2{x}+b^2\sin^2{x}}}=
\frac{1}{a}\int_0^{\phi}\frac{\mathrm{d}x}{\sqrt{1-\left(1-\frac{b^2}{a^2}\right)\sin^2{x}}}=
\end{equation}
$$
=\frac{1}{a}F\left(\phi,1-\frac{b^2}{a^2}\right)
$$
and
$$
I(a,b,\frac{\pi}{2})=
\frac{1}{a}K\left(1-\frac{b^2}{a^2}\right).
$$
Thus, the equality (\ref{eq1}) may be rewritten as $\frac{1}{M(a,b)}=\frac{2}{a\pi}K\left(1-\frac{b^2}{a^2}\right)$
or
$$
K\left(1-\frac{b^2}{a^2}\right)=\frac{\pi}{2}\frac{1}{\frac{1}{a}M(a,b)}=\frac{\pi}{2}\frac{1}{M(1,\frac{b}{a})}.
$$

As in \cite{1},  for $I(a,b,\phi)$  the changing of variables
$$
\sin{x}=\frac{2a\sin{\varphi}}{a+b+(a-b)\sin^2{\varphi}}
$$
leads to the sequence
\begin{equation}\label{eq11}
I(a,b,\phi)\stackrel{def}{=}I_0(a_0,b_0,\phi_0)=I_1(a_1,b_1,\phi_1)=I_2(a_2,b_2,\phi_2)=\ldots
\end{equation}
where
$$
I_k(a_k,b_k,\phi_k)=\int_0^{\phi_k}\frac{\mathrm{d}\varphi}{\sqrt{a_k^2\cos^2{\varphi}+b_k^2\sin^2{\varphi}}}
$$
and the upper integration limits are  generated by the sequence
$$
\sin{\phi_{k-1}}=\frac{2a_{k-1}\sin{\phi_k}}{a_{k-1}+b_{k-1}+(a_{k-1}-b_{k-1})\sin^2{\phi_k}}.
$$
The sequence $(\sin{\phi_k})_{k\in\mathbb{N}}$ is decreasing and consequently
the sequence $(\phi_k)_{k\in\mathbb{N}}$ is convergent.
It results that
\begin{eqnarray}
\sin{\phi_k}&=&a_{k-1}-\frac{\sqrt{a_{k-1}^2\cos^2{\phi_{k-1}}+b_{k-1}^2\sin^2{\phi_{k-1}}}}
{(a_{k-1}-b_{k-1})\sin{\phi_{k-1}}}=y_k \label{eq14}\\
\phi_k &=&\arcsin{y_k}. \nonumber
\end{eqnarray}
From (\ref{eq11}) it results 
$$
I(a,b,\phi)=\lim_{k\rightarrow\infty}I_k(a_k,b_k,\phi_k)=\frac{\phi_{\infty}}{M(a,b)},
$$
with $\phi_{\infty}=\lim_{k\rightarrow\infty}\phi_k.$ Using (\ref{eq1}) we get
$$
I(a,b,\phi)=\frac{1}{a}F\left(\phi,1-\frac{b^2}{a^2}\right)=\frac{\phi_{\infty}}{M(a,b)}
$$
and consequently
$$
F\left(\phi,1-\frac{b^2}{a^2}\right)=\frac{a\phi_{\infty}}{M(a,b)}=\frac{\phi_{\infty}}{\frac{1}{a}M(a,b)}
=\frac{\phi_{\infty}}{M(1,\frac{b}{a})}.
$$
Denoting $m=1-\frac{b^2}{a^2},  (a>b>0\ \Leftrightarrow\ 0<m<1),$ the above equation becomes
\begin{equation}\label{eq31}
F(\phi,m)=\frac{\phi_{\infty}}{M(1,\sqrt{1-m})}.
\end{equation}
Therefore the computation of $F(\phi,m)$ returns to generate iteratively the sequences $(a_k)_k,\ (b_k)_k,\
(\phi_k)_k$ until a stopping condition is fulfilled. The initial values are $a_0=1,\ b_0=\sqrt{1-m}$ and 
$\phi_0=\phi.$ For $a_0=1,$ instead of the sequences $(a_k)_k,\ (b_k)_k$ we may compute the sequences, \cite{6},
$$
\begin{array}{l}
s_0 =b_0\\
s_{k+1} = \frac{2\sqrt{s_k}}{1+s_k}
\end{array}
\qquad
\begin{array}{l}
p_0 =\frac{1}{2}(1+s_0)\\
p_{k+1} = \frac{1}{2}(1+s_k) p_k
\end{array}.
$$
Then $\lim_{k\rightarrow\infty}p_k=M(1,b_0).$

If $\phi=\frac{\pi}{2}$ then $\phi_{\infty}=\frac{\pi}{2}$ and we retrieve 
\begin{equation}\label{eq30} 
K(m)=\frac{\pi}{2M(1,\sqrt{1-m})}.
\end{equation}

From a practical point of view and as a drawback the method is not applicable when $\phi$ is small, e.g. $0<\phi<10^{-5}.$
The cause is the presence of the factor $\sin{\phi_{k-1}}$ in the denominator in (\ref{eq14}).
In this case, from the Maclaurin series expansion of $F(\phi,m)$  we get
$F(\phi,m)\approx \phi-\frac{m}{6}\phi^3.$

\section{The Jacobi elliptic function \texttt{sn}}

The Jacobi elliptic function $\sn(z,m)$ may be defined by the equation, \cite{5},
\begin{equation}\label{e12}
z=\int_0^{\sn(z,m)}\frac{\mathrm{d}t}{\sqrt{(1-t^2)(1-mt^2)}}.
\end{equation}

Throughout this paper the variable $m$ is fixed and we use the shorter notation $\sn(z),$ omitting $m.$

We shall use the following Jacobi elliptic functions, too
$$
\cn^2(z)=1-\sn^2(z), \qquad \dn^2(z,m)=1-m\ \sn^2(z).
$$
Again we shall use the shorter notation $\dn(z).$

If
$$
K=K(m) \qquad \mbox{and}\qquad K'=K(1-m)
$$ 
then the parallelogram period is the rectangle $D=[0,4K)+i[0,2K')$  and the points
$K'i$ and $2K+iK'$ are poles of first order, \cite{6}.

The following properties of the function $\sn$ will be used, \cite{6}, \cite{7}:
\begin{itemize}
\item
\begin{equation}\label{eq15}
\sn(-z)=-\sn(z)
\end{equation}
\item
\begin{eqnarray}
\sn(x+y)&=&\frac{\sn(x)\cn(y)\dn(y)+\sn(y)\cn(x)\dn(x)}{1-m\ \sn^2(x)\sn^2(y)}\label{eq16}\\
\cn(x+y)&=&\frac{\cn(x)\cn(y)-\sn(x)\sn(y)\dn(x)\dn(y)}{1-m\ \sn^2(x)\sn^2(y)}\label{eq21}\\
\dn(x+y)&=&\frac{\dn(x)\dn(y)-m\ \sn(x)\sn(y)\cn(x)\cn(y)}{1-m\ \sn^2(x)\sn^2(y)}\label{eq22}
\end{eqnarray}
\item
Because $\sn(K)=1,\ \cn(K)=0,\ \dn(K)=\sqrt{1-m}$ from the above equalities it results
\begin{eqnarray}
\sn(K\pm z)&=&\frac{\cn(z)}{\dn(z)}\label{eq17}\\
\cn(K+z)&=&-\sqrt{1-m}\frac{\sn(z)}{\dn(z)}\label{eq23}\\
\dn(K+z)&=&\frac{\sqrt{1-m}}{\dn(z)}\label{eq24}
\end{eqnarray}
\item
Knowing that $\sn(2K)=0, \cn(2K)=-1, \dn(2K)=1$, from (\ref{eq16}) it results
\begin{equation}\label{eq18}
\sn(2K\pm z)=\mp\sn(z)
\end{equation}
\item
Knowing that $\sn(K+iK')=\frac{1}{\sqrt{m}}, \dn(K+iK')=0$ from (\ref{eq16}) it results
\begin{equation}\label{eq20}
\sn(z+K+iK')=\frac{1}{\sqrt{m}}\ \frac{\dn(z)}{\cn(z)}.
\end{equation}
\end{itemize}

The computation of $\sn(z)$ depends on the position of $z$ in $D$ and we suppose
that we know $K$ and $K'.$
\begin{itemize}
\item
If $z\in [0,K]$ or $z\in i[0,K')$ then $\sn(z)$ will be the solution $u$ of the equation
\begin{equation}\label{eq13}
\Phi(u)=\int_0^u\frac{\mathrm{d}t}{\sqrt{(1-t^2)(1-mt^2)}}-z=0.
\end{equation}
\item
Otherwise and excepting the poles the value of $\cn(z)$ will be computed
using the properties of the function $\sn$ and
its values on some points from $[0,K]\cup i[0,K'),$ without any other additional
numerical procedure. 
\end{itemize}

\subsection*{Computing $\sn(z)$ in the segment $[0,4K)$}
The following cases arise:
\begin{enumerate}
\item
$z\in [0,K].$ Equation (\ref{eq13}) may be solved using the Newton method with the iterations
$$
u_{k+1}=u_k-\frac{\Phi(u_k)}{\Phi'(u_k)}=u_k-\left(F(\arcsin{u_k},m)-z\right)\sqrt{(1-u_k^2)(1-mu_k^2)}
$$
The linear interpolation between $\sn(z,0)$ and $\sn(z,1)$ gives the initial approximation 
$u_0=(1-m)\sin{z}+m\tanh{z}.$

If $z$ is small enough the method is rapidly converging and for $z$ near $K$ we set $z'=K-z$ and
after computing $\sn(z')=w'$ as was described above, using (\ref{eq17}) we have  
$$
\sn(z)=\sn(K-z')=\frac{\cn(z')}{\dn(z')}=\sqrt{\frac{1-w'^2}{1-m\ w'^2}}.
$$
\item
$z\in (K,4K).$
Let be
$$
z'=\left\{\begin{array}{lcl}
2K-z & \mbox{if} & z\in (K,2K]\\
z-2K &  \mbox{if} & z\in (2K,3K]\\
4K-z &  \mbox{if} & z\in (3K,4K)
\end{array}\right..
$$
After computing $\sn(z')=w',\ z'\in [0,K],$ we have
$$
\sn(z)=\left\{\begin{array}{lcl}
w' & \mbox{if} & z\in (K,2K]\\
-w' &  \mbox{if} & z\in (2K,4K)
\end{array}\right..
$$
Indeed, if $z\in (K,2K]$ then 
$$
\sn(z)=\sn(2K-z')=\sn(z')=w';
$$
if $z\in(2K,3K]$ then
$$
\sn(z)=\sn(2K+z')=-\sn(z')=-w'
$$
and if $z\in(3K,4K)$ then
$$
\sn(z)=\sn(4K-z')=\sn(-z')=-\sn(z')=-w'.
$$
\end{enumerate}

\subsection*{Computing $\sn(z)$ for $z\in i(0,2K')\setminus\{iK'\}.$}
The following cases arise:
\begin{enumerate}
\item
$z\in i[0,K').$ Writing $z=iy,\ y\in [0,K'),$ we are looking for the solution of the equation (\ref{eq13}) in the form
$u=iv, v\in\mathbb{R}. $ After the change of variable $t=i s$ there is obtained the equation
\begin{equation}\label{eq19}
\Psi(v)=\int_0^v\frac{\mathrm{d}s}{\sqrt{(1+s^2)(1+m\ s^2)}}-y=0.
\end{equation}
According to the Newton method, the iterations are
$$
v_{k+1}=v_k-\frac{\Psi(v_k)}{\Psi'(v_k)}=
$$
$$
=v_k-\left(\int_0^{v_k}\frac{\mathrm{d}s}{\sqrt{(1+s^2)(1+m\ s^2)}}-y\right)
\sqrt{(1+v_k^2)(1+m\ v_k^2)}
$$ 
starting with $v_0=y.$ The above integral is computed using a quatrature procedure.
\item
$z\in i(K',2K'].$ Let be $z'=z-iK'=iy'$ with $y'\in(0,K').$ From (\ref{eq20}) we have
$$
\sn(z)=\sn(iK'+iy')=\frac{1}{\sqrt{m}}\ \frac{\dn(iy'-K)}{\cn(iy'-K)}.
$$
After using (\ref{eq23}) and (\ref{eq24}) it results
$$
\sn(z)=\frac{1}{\sqrt{m}\ \sn(z')}.
$$
\end{enumerate}

\subsection*{Computing $\sn(z)$ in the rectangle period}

We describe here how to compute $\sn(z)$ when $z$ belongs to the rectangle period
excepting the poles and the lower and the left sides.

Let $z=x+iy$ such that
$x\in(0,4K)$ and $y\in (0,2K').$

The following cases arise:
\begin{enumerate}
\item
$y\not=K'.$ Using (\ref{eq16}) we have
$$
\sn(z)=\sn(x+iy)=\frac{\sn(x)\cn(iy)\dn(iy)+\sn(iy)\cn(x)\dn(x)}{1-m\ \sn^2(x)\sn^2(iy)}.
$$
$\sn(x), \sn(iy)$ are computed as was presented above and then compute
$\cn(x), \cn(iy), \dn(x), \dn(iy).$ It must be taken into account that if $x\in(K,3K)$ then
$\cn(x)=-\sqrt{1-\sn^2(x)}.$
\item
$y=K'.$ We deduce through (\ref{eq24})
$$
\sn(z)=\sn(K+iK'+x-K)=\frac{1}{\sqrt{m}}\ \frac{\dn(x-K)}{\cn(x-K)}.
$$
As above, applying (\ref{eq23}) and (\ref{eq24}) we obtain
$$
\sn(z)=\frac{1}{\sqrt{m}\ \sn(x)}.
$$
\end{enumerate}
On poles, $z\in\{iK',2K+iK'\},$ we set  $\sn(z)=\infty.$

In his way we have computed the value of $\sn(z)$ for any $z\in D.$

\section{Numerical results}

Numeric computing softwares contains functions to compute elliptical integrals and elliptical functions.
We recall some methods in \textit{Mathematica} and \textit{Scilab} in Table 1.

\begin{table}[h]
\hspace*{4cm}
\begin{tabular}{|l|l|}
\hline\hline
  Meaning  & Function signature\\
\hline\hline
\multicolumn{2}{|c|}{Mathematica}\\
\hline\hline
$F(\phi,m)$ & EllipticF[$\phi,m$] \\
\hline
$K(m)$ & EllipticK[$m$] \\
\hline
$\sn(z,m)$ & JacobiSN[$z,m$] \\
\hline
$\cn(z,m)$ & JacobiCN[$z,m$] \\
\hline
$\dn(z,m)$ & JacobiDN[$z,m$] \\
\hline\hline
\multicolumn{2}{|c|}{Scilab}\\
\hline\hline
$F(\phi,m)$ & delip($\phi,\sqrt{m}$)\\
\hline
$\sn(z,m)$ & \%$\sn(z,m)$\\
\hline\hline
\end{tabular}

\caption{Elliptical function in \textit{Mathematica} and \textit{Scilab}.}
\end{table}

We developed a \textit{Scilab} program based on the method presented in this paper. 
The values $K$ and $K'$ were computed using (\ref{eq30}). 
Because the numbers have a floating point representation two numbers are considered to be equal 
if their distance is less than a tolerance.

Some results are given in Table 2. The function \texttt{JacobiSN} from
\textit{Mathematica} gives similar values (excepting the poles).

\begin{landscape}
\begin{table}[h]
\hspace*{4cm}
\begin{tabular}{|l|l|l|l|}
\hline\hline
$z$ & $\sn(z)$ computed &  \%$\sn(z,m)$ & Error \\
 & & & $|\sn(z)-\%\sn(z,m)|$ \\
\hline\hline
\multicolumn{4}{|c|}{$m=0.81$} \\
\hline
K=2.2805491  & \multicolumn{2}{|l|}{delip$(1,0.9)=2.2805491$} & $4.441D-16$ \\
\hline
K'=1.6546167 & \multicolumn{2}{|l|}{delip$(1,\sqrt{0.19})=1.6546167$} & $2.220D-16$\\
\hline\hline
$0.5K$  &   0.8345252    &     0.8345252       &         1.024D-09    \\
$1.4K$  &   0.9038225    &  0.9038225         &       6.602D-10\\
$2.7K$  &   -0.9501563    &     -0.9501563      &          7.363D-10 \\
$3.3K$   &  -0.9501563    &     -0.9501563      &          7.363D-10   \\
$i 0.6K'$   &   1.4511449i   &    1.4511449i      &         1.554D-15   \\
$i 1.3K' $    &    -2.0696167i   &  -2.0696167i     &          1.332D-15  \\     
$0.8K+i 0.3K'$  &  1.0085488 + 0.0420829i & 1.0085488 + 0.0420829i &  5.812D-10\\
$0.5K+i 1.7K'$ &  0.9048397 - 0.1679796i  & 0.9048397 - 0.1679796i  & 1.129D-09\\
$1.3K+i 1.7K'$ &  0.9892195 + 0.071665i  & 0.9892195 + 0.071665i &   2.701D-10\\
$2.5K+i 0.4K'$ & -0.9592212 - 0.2093038i  & -0.9592212 - 0.2093038i &  8.724D-10\\
$3.6K+i 0.4K;$ & -0.8951883 + 0.3091877i  &  -0.8951883 + 0.3091877i &  1.656D-09\\
$3.6K+i 1.7K'$ & -0.8233279 - 0.2419397i  & -0.8233279 - 0.2419397i &  1.519D-09\\
$0.5K+i K'$ &  1.3314291  &    1.3314291       &         1.634D-09       \\
$2.5K+i K'$ &  -1.3314291 &     -1.3314291     &           1.183D-09       \\
$K+i K'$ &  1.1111111      &    1.1111111      &          2.220D-16    \\
$K$       &         1.       &     1.          &             2.220D-16         \\
$i K'$     &         Nan + Infi  &       1.633D+16i       &        Nan    \\      
$2K+i K'$ &  Nan + Infi   &    -4.211D+15 + 1.170D+15i  &  Nan    \\
\hline\hline    
\end{tabular}

\caption{Results obtained using the presented method.}
\end{table}
\end{landscape}

The 3D image of the modules of the function $\sn(z,0.81)$ computed on $D$ is given in Figure 1.

Finally we show the visualization of the complex function, using the method presented in  \cite{9}. In a point
the value of the function is represented by a color obtained \textit{projecting}  that value
into the colors cube. The procedure is based on the stereographic projection. 

The Figure 2 is given for calibration, representing the visualization of the identity function. 

The Figure 3
contains the visualization of $\sn(z,0.81),\ z\in D.$ The zeros are colored in black while
the poles are colored in white.

\begin{center}
\begin{figure}[h]
\hspace*{2cm}\includegraphics[width=13cm,height=6cm,keepaspectratio]{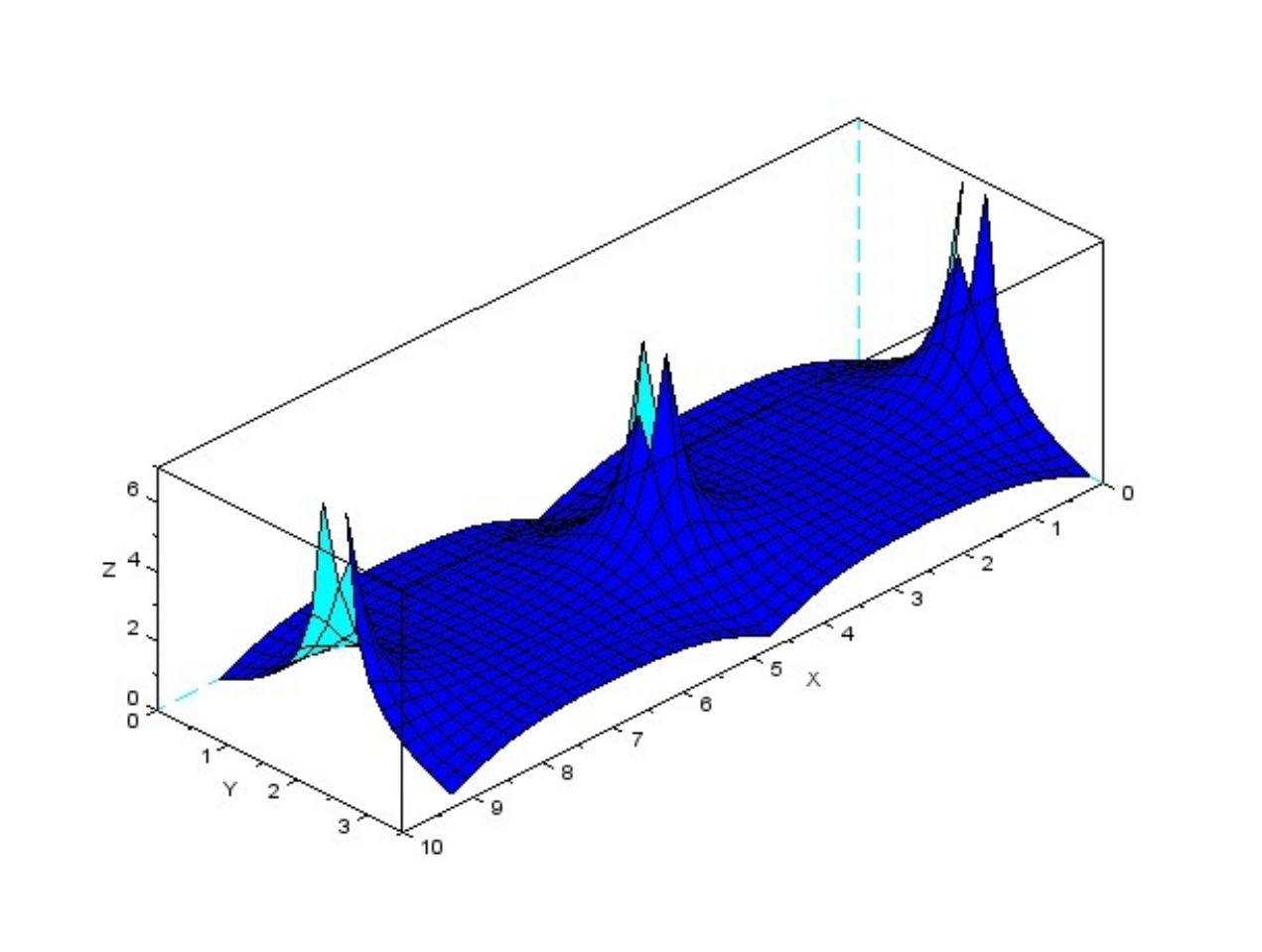}

\caption{For $m=0.81$ the 3D image of the modules of $\sn$.}
\end{figure}
\end{center}

\begin{figure}[h]
\includegraphics[width=13cm,height=6cm,keepaspectratio]{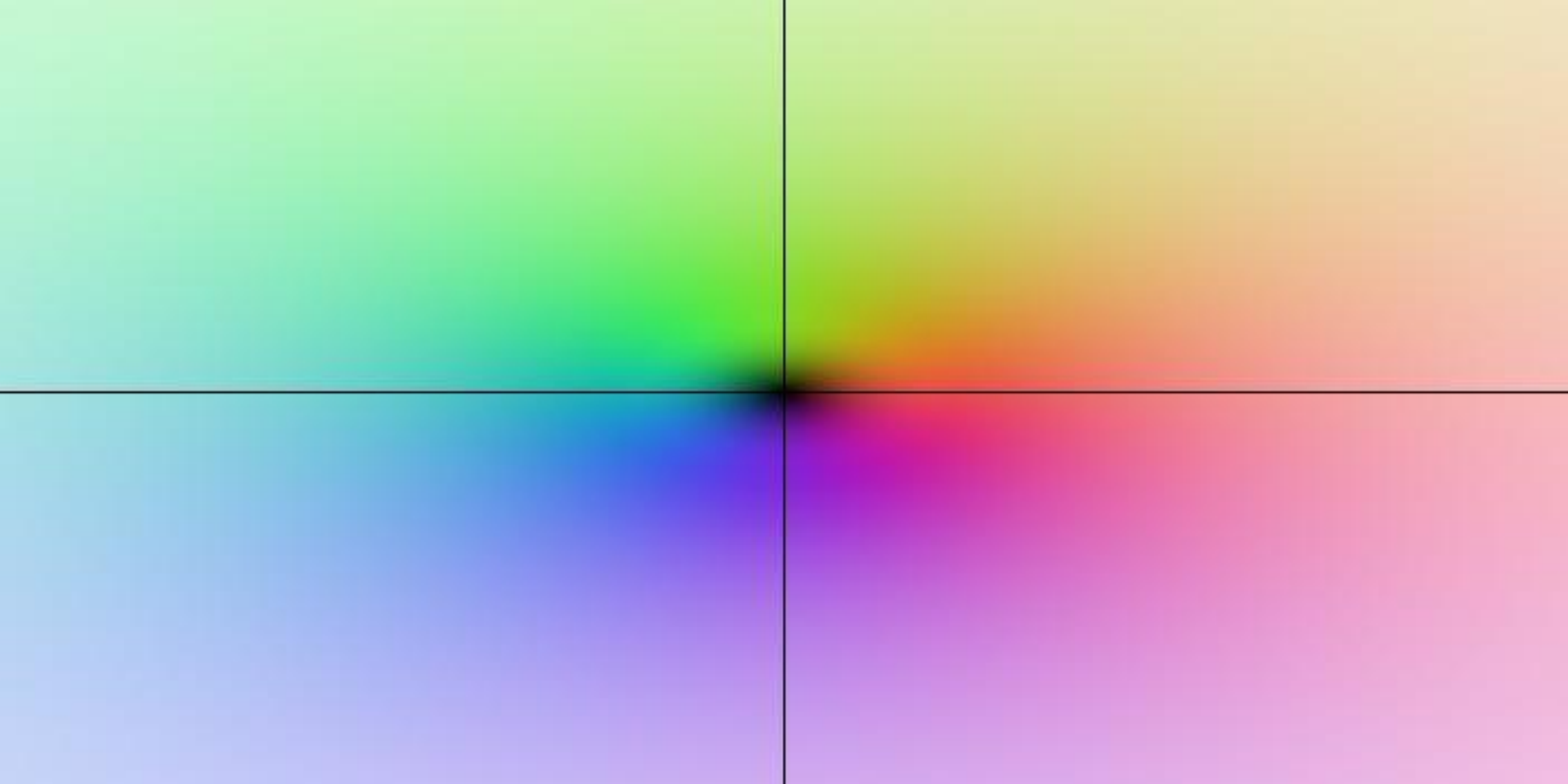}

\caption{Function $w\mapsto w, w\in [-5,5]+i[-5,5]$}

\vspace*{0.3cm}
\includegraphics[width=13cm,height=6cm,keepaspectratio]{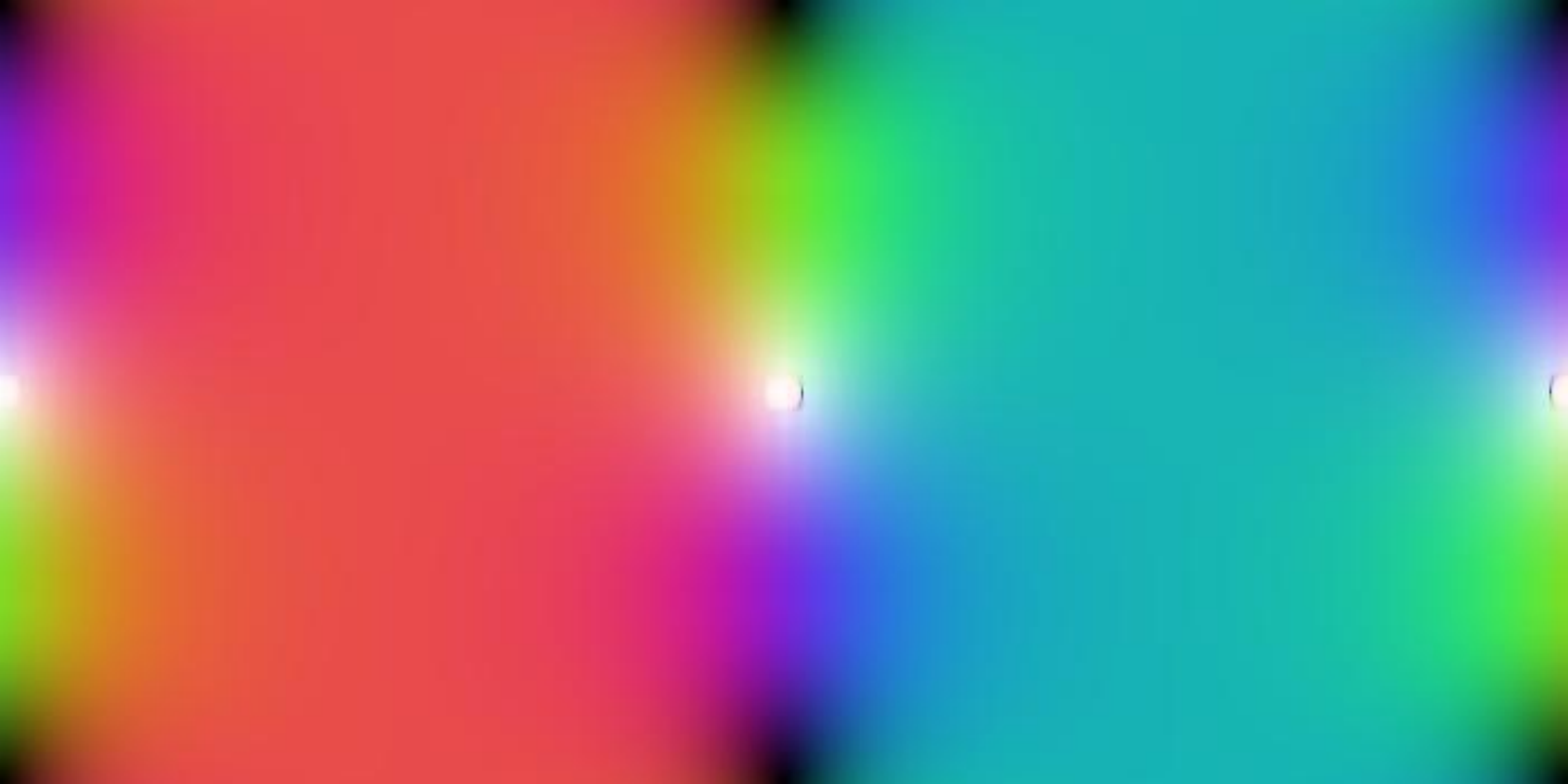}

\caption{Function $w=\sn(z,0.81), z\in D$.}
\end{figure}

\end{document}